\documentclass[12pt]{article}
\usepackage{amsthm}
\usepackage{mathrsfs}
\usepackage{graphicx}
\usepackage{amsmath}
\usepackage{amsfonts}
\usepackage{amssymb}
\usepackage{url}
\usepackage{hyperref}
\usepackage{dsfont} 
\usepackage{verbatim}
\usepackage[mathlines]{lineno}
\usepackage{color}
\usepackage{enumerate}
\usepackage{enumitem}
\usepackage{graphicx}
\definecolor{red}{rgb}{1,0,0}

\definecolor{pink}{rgb}{.9,.2,.7}

\definecolor{org}{rgb}{1,.7,0}

\definecolor{purp}{rgb}{.5,0,.5}

\definecolor{br}{rgb}{.5,.7,1}

\newtheorem{thm}{Theorem}[section]

\newtheorem{prop}[thm]{Proposition}
\newtheorem{cor}[thm]{Corollary}
\newtheorem{rem}[thm]{Remark}
\newtheorem{lem}[thm]{Lemma}
\newtheorem{ex}[thm]{Example}

\newtheorem{obs}[thm]{Observation}

\def\noi{\noindent}
\def\mtx#1{\begin{bmatrix} #1 \end{bmatrix}}

\newcommand{\Z}{\mathbb{Z}}

\newcommand{\ba}{\begin{array}}
\newcommand{\ea}{\end{array}}

\newcommand{\bit}{\begin{itemize}}
\newcommand{\eit}{\end{itemize}}
\newcommand{\ben}{\begin{enumerate}}
\newcommand{\een}{\end{enumerate}}
\newcommand{\bea}{\begin{eqnarray*}}
\newcommand{\eea}{\end{eqnarray*}}
\newcommand{\beq}{\begin{equation}}
\newcommand{\eeq}{\end{equation}}
\newcommand{\bpf}{\begin{proof}}
\newcommand{\epf}{\end{proof}}

\DeclareMathSymbol{\mlq}{\mathord}{operators}{``}
\DeclareMathSymbol{\mrq}{\mathord}{operators}{`'}

\newcommand{\OL}{\overline}
\newcommand{\rank}{\operatorname{rank}}

\newcommand{\epr}{\operatorname{epr}}
\newcommand{\pr}{\operatorname{pr}}
\newcommand{\ones}{\ensuremath{\mathds{1}}}

\newcommand{\F}{\mathbb{F}}

\begin{document}


\title{The enhanced principal rank characteristic sequence over a field of characteristic $2$}
\author{Xavier Mart\'{i}nez-Rivera\thanks{Department of Mathematics, Iowa State University, Ames,
IA 50011, USA (xaviermr@iastate.edu).}
}

\maketitle

\begin{abstract}
The enhanced principal rank characteristic sequence
(epr-sequence) of an $n \times n$ symmetric matrix
over a field $\F$ was recently defined as
$\ell_1 \ell_2 \cdots \ell_n$, where $\ell_k$ is either
$\tt A$, $\tt S$, or $\tt N$ based on whether
all, some (but not all),  or none of the
order-$k$ principal minors of the matrix are nonzero.
Here, a complete characterization of the epr-sequences that are attainable by symmetric matrices over the field $\Z_2$, the integers modulo $2$, is established.
Contrary to the attainable epr-sequences over a field of characteristic $0$, our characterization reveals that the attainable epr-sequences over $\Z_2$ possess very special structures.
For more general fields of characteristic $2$,
some restrictions on attainable epr-sequences are obtained.
\end{abstract}

\noi{\bf Keywords.}
Principal rank characteristic sequence;
enhanced principal rank characteristic sequence;
minor; rank; symmetric matrix; finite field.\\

\noi{\bf AMS subject classifications.}
15A15, 15A03.


\section{Introduction}
$\null$
\indent
For an $n \times n$ real symmetric matrix $B$,
Brualdi et al.\ \cite{ORIGINAL} introduced the
\textit{principal rank characteristic sequence}
(abbreviated pr-sequence), which was defined as
$\pr(B) = r_0]r_1 \cdots r_n$, where, for $k \geq 1$,
   \begin{equation*}
      r_k =
         \begin{cases}
             1 &\text{if $B$ has a nonzero principal minor of order $k$, and}\\
             0 &\text{otherwise,}
         \end{cases}
   \end{equation*}
while $r_0 = 1$ if and only if $B$ has a
$0$ diagonal entry.
This definition was generalized for symmetric matrices
over any field by Barrett et al.\ \cite{FIELDS}.

Our focus will be studying a sequence
that was introduced by Butler et al.\ \cite{EPR} as a refinement of the pr-sequence of an $n \times n$
symmetric matrix $B$ over a field $\F$,
which they called the
\textit{enhanced principal rank characteristic sequence}
(abbreviated epr-sequence),
and which was defined as
$\epr(B) = \ell_1\ell_2 \cdots \ell_n$, where
   \begin{equation*}
      \ell_k =
         \begin{cases}
             \tt{A} &\text{if all the principal minors of order $k$ are nonzero;}\\
             \tt{S} &\text{if some but not all the principal minors of order $k$ are nonzero;}\\
             \tt{N} &\text{if none of the principal minors of order $k$ are nonzero, i.e., all are zero.}
         \end{cases}
   \end{equation*}
The definition of the epr-sequence was later extended to the class of real skew-symmetric matrices in \cite{skew}, where a complete characterization of the epr-sequences realized by this class was presented.
However, things are more subtle for the class of symmetric matrices over a field $\F$, and thus obtaining a similar characterization presents a difficult problem.
When $\F$ is of characteristic $0$,
it is known that any epr-sequence of the form
$\ell_1 \cdots \ell_{n-k} \tt \OL{N}$,
with $\ell_i \in \{\tt A,S\}$,
is attainable by an $n \times n$ symmetric matrix over $\F$, where $\OL{\tt N}$ (which may be empty) is the sequence consisting of
$k$ consecutive  $\tt N$s \cite{EPR} -- if $\OL{\tt N}$ is empty, note that we must have $\ell_n = \tt A$.
In general, the subtlety for symmetric matrices becomes evident once the $\tt N$s are not restricted to occur consecutively at the end of the sequence:
Sequences such as $\tt NSA$, $\tt NNA$ and $\tt NNS$ can never occur as a subsequence of the epr-sequence of a symmetric matrix over any field \cite{EPR};
the same holds for the sequences $\tt NAN$ and $\tt NAS$ when the field is of characteristic not $2$ \cite{EPR}.
Moreover, over fields of characteristic not $2$,
the sequence $\tt ANS$ can only occur at the start of the sequence \cite{EPR}.
Over the real field, $\tt SNA$ can only occur as
a terminal subsequence, or in
the terminal subsequence $\tt SNAA$ \cite{XMR-Classif}.
Furthermore, over the real field, we also know that when the subsequence $\tt ANA$ occurs as a non-terminal subsequence, it forces every other term of the sequence to be $\tt A$ \cite{XMR-Classif}.
However, it is unknown what kind of restrictions a subsequence such as $\tt SNS$ imposes on an attainable sequence (over any field);
this is one of the difficulties in arriving at a complete characterization of the epr-sequences attainable by a symmetric matrix over a field $\F$.
In order to simplify this problem,
it is natural to consider the case when
$\F$ is of characteristic $2$.
The analogous problem for pr-sequences was
already settled in \cite{FIELDS}:

\begin{thm}\label{pr-char2 classif}
{\rm \cite[Theorem 3.1]{FIELDS}}
A pr-sequence of order $n \geq 2$ is attainable by a
symmetric matrix over a field of characteristic $2$
if and only if
it has one of the following forms:
\[
0]1 \hspace{0.1cm} \OL{1} \hspace{0.1cm} \OL{0},
\qquad
1]  \OL{01} \hspace{0.1cm} \OL{0},
\qquad
1]1 \hspace{0.1cm} \OL{1} \hspace{0.1cm} \OL{0}.
\]
\end{thm}

We see that for any two fields of characteristic $2$, the class of pr-sequences attainable by symmetric matrices over each of the two fields is the same.
This is not true in the case of epr-sequences:
Consider an epr-sequence starting with $\tt AA$ over the field $\Z_2 = \{0,1\}$, the integers modulo $2$;
over this field, any such sequence must be $\tt AA\OL{A}$, since any symmetric matrix attaining this sequence must be the identity matrix.
However, in Example \ref{ex: counterexamples} below, it is shown that the epr-sequence $\tt AAN$ is attainable over a field of characteristic $2$, implying that not all fields of characteristic $2$ give rise to the same class of attainable epr-sequences.
In light of this difficulty,
our main focus here will be on the field $\F = \Z_2$;
after establishing some restrictions for
the attainability of epr-sequences over a field of characteristic $2$ at the beginning of
Section \ref{s: Restrictions over any field of char 2},
our main objective is a complete
characterization of the epr-sequences that
are attainable by symmetric matrices over $\Z_2$
(see Theorems \ref{N...}, \ref{A...} and \ref{S...}).
We find that the attainable epr-sequences over $\Z_2$ possess very special structures,
which is in contrast to the family of
attainable epr-sequences over a field of
characteristic $0$, which was described above.

Another motivating factor for considering this problem is
that it is a simplification of the
\textit{principal minor assignment problem}
as stated in \cite{PMAP introduced},
which also served as motivation for the introduction
of the pr-sequence in \cite{ORIGINAL}.
Note that epr-sequences provide more information than
pr-sequences, and thus are a step closer to
the principal minor assignment problem.

Extra motivation for this problem comes from the observation that there is a one-to-one correspondence between adjacency matrices of simple graphs and symmetric matrices over $\Z_2$ with zero diagonal,
and, more generally, between adjacency matrices of loop graphs and symmetric matrices over $\Z_2$.

It should be noted that,
although epr-sequences have received attention after
their introduction in \cite{EPR}
(see
\cite{EPR-Hermitian},
\cite{skew} and
\cite{XMR-Classif}, for example),
very little is known about epr-sequences of
symmetric matrices over a field of
characteristic $2$, since the vast majority of
what has appeared on the literature regarding
epr-sequences has been focused on fields of characteristic
\textit{not} $2$.

Although Theorem \ref{pr-char2 classif} sheds some
light towards settling the problem under consideration,
it does not render it trivial by any means;
one reason is the observation that
two symmetric matrices may have distinct epr-sequences
while having the same pr-sequence:
As it is shown in Theorem \ref{A...} below,
the epr-sequences $\tt ASAA$ and $\tt ASSA$,
which are associated with the pr-sequence $0]1111$,
are both attainable over $\Z_2$.

To highlight a second reason,
we state the two results upon which
Barrett et al.\ \cite{FIELDS} relied in order to obtain
Theorem \ref{pr-char2 classif}
(the latter is a variation of a result of
Friedland \cite[p. 426]{Friedland}).

\begin{lem}\label{pr-char2 lemma 1}
{\rm \cite[Lemma 3.2]{FIELDS}}
Let $\F$ be a field of characteristic $2$,
let $B$ be a symmetric matrix over $\F$
with $\pr(B) = r_0] r_1 \cdots r_n$, and
let $E$ be an $n \times n$ invertible matrix over $\F$.
Then $\epr(EBE^T) = r'_0] r_1 r_2 \cdots r_n$
for some $r'_0 \in \{0,1\}$.
\end{lem}

In what follows, $K_n$ denotes
the complete graph on $n$ vertices, and
$A(K_n)$ denotes its adjacency matrix.

\begin{lem}\label{pr-char2 lemma 2}
{\rm \cite[Lemma 3.3]{FIELDS}}
Let $B$ be a symmetric matrix over a field $\F$
with characteristic $2$.
Then $B$ is congruent to the direct sum of a
(possibly empty) invertible diagonal matrix $D$,
and a (possibly empty) direct sum of $A(K_2)$ matrices,
and a (possibly empty) zero matrix.
\end{lem}

The two lemmas above permitted
Barrett et al.\ \cite{FIELDS} to arrive at
their characterization for pr-sequences
in Theorem \ref{pr-char2 classif} by
restricting themselves to symmetric matrices that
are in the canonical form described in
Lemma \ref{pr-char2 lemma 2}.
We cannot use this approach to obtain our desired  characterization for epr-sequences:
Suppose one tries to apply the congruence described in
Lemma \ref{pr-char2 lemma 1} to a symmetric matrix $B$
with $\epr(B) = \tt ASAN$,
which is shown to be attainable in Theorem \ref{A...}.
Then, because $B$ is singular,
and because multiplication by an invertible matrix
preserves the rank of the original matrix,
once $B$ has been transformed into the
canonical form described in Lemma \ref{pr-char2 lemma 2},
it must be the case that in this resulting matrix
the zero summand is non-empty.
Thus, the resulting matrix has a zero row
(and zero column),
which implies that it contains a principal minor of
order $3$ that is zero.
Then, as the principal minors of order $3$ of the
original matrix $B$ were all nonzero,
the congruence performed did not preserve the
third term of $\epr(B)$,
which is in contrast to what happens to $\pr(B)$,
which, with the exception of the zeroth term,
must be preserved completely by
Lemma \ref{pr-char2 lemma 1}.

We say that a (pr- or epr-) sequence is
\textit{attainable} over a field $\F$ provided that there exists a symmetric matrix $B \in \F^{n \times n}$ that attains it.
A pr-sequence and an epr-sequence are
\textit{associated} with each other if a matrix
(which may not exist) attaining the epr-sequence also
attains the pr-sequence.
A subsequence that does not appear in any attainable
sequence is \textit{prohibited}.
We say that a sequence has \textit{order} $n$ if it corresponds to a matrix of order $n$.
Let $B$ be an $n \times n$ matrix, and let
$\alpha, \beta \subseteq \{1, 2, \dots, n\}$;
then the submatrix lying in rows indexed by $\alpha$,
and columns indexed by $\beta$,
is denoted by $B[\alpha, \beta]$.
The matrix obtained by deleting the rows indexed by
$\alpha$, and columns indexed by $\beta$,
is denoted by $B(\alpha, \beta)$.
If $\alpha = \beta$, then the principal submatrix
$B[\alpha, \alpha]$ is abbreviated to $B[\alpha]$,
while $B(\alpha, \alpha)$ is abbreviated to $B(\alpha)$.
The matrices $O_n$ and $I_n$ denote, respectively,
the zero and identity matrix of order $n$.
We denote by $J_{m,n}$ the $m \times n$ all-$1$s matrix,
and, when $m=n$, $J_{n,n}$ is abbreviated to $J_n$.
The block diagonal matrix formed from
two square matrices $B$ and $C$ is
denoted by $B \oplus C$.
The matrices $B$ and $C$ are
\emph{permutationally} similar if there exists a permutation matrix $P$ such that $C=P^TBP$.
Given a graph \textit{G}, $A(G)$ denotes the adjacency matrix of $G$.

\subsection{Results cited}
$\null$
\indent
This section lists results that
will be cited frequently, with some of them being assigned abbreviated nomenclature.

\begin{thm}\label{NN result}
{\rm \cite[Theorem 2.3]{EPR} ($\tt NN$ Theorem.)}
Suppose $B$ is a symmetric matrix over
a field $\F$,
$\epr(B) = \ell_1\ell_2 \cdots \ell_n$, and
$\ell_k = \ell_{k+1} = \tt{N}$ for some $k$.
Then $\ell_i = \tt{N}$ for all $i \geq k$.
\end{thm}

\begin{thm}\label{Inverse Thm}
{\rm\cite[Theorem 2.4]{EPR} (Inverse Theorem.)}
Suppose $B$ is a nonsingular symmetric matrix
over a field $\F$.
If $\epr(B) = \ell_1\ell_2 \cdots \ell_{n-1}\tt{A}$, then
$\epr(B^{-1}) = \ell_{n-1}\ell_{n-2} \cdots \ell_{1}\tt{A}$.
\end{thm}

Given a matrix $B$, the $i$th term in its epr-sequence is
denoted by $[\epr(B)]_i$.

\begin{thm}\label{Inheritance}
{\rm{\cite[Theorem 2.6]{EPR}}}
{\rm (Inheritance Theorem.)}
Suppose that  $B$ is a symmetric matrix over a field $\F$, $m \leq n$, and $1\le i \le m$.  
\ben
\item
If $[\epr(B)]_i={\tt N}$, then  $[\epr(C)]_i={\tt N}$ for all $m\times m$ principal submatrices $C$.
\item If  $[\epr(B)]_i={\tt A}$, then  $[\epr(C)]_i={\tt A}$ for all $m\times m$ principal submatrices $C$.

\item If $[\epr(B)]_m={\tt S}$, then there exist $m\times m$ principal submatrices $C_A$ and $C_N$ of $B$ such that $[\epr(C_A)]_m = {\tt A}$ and $[\epr(C_N)]_m = {\tt N}$.
\item If $i < m$ and $[\epr(B)]_i = {\tt S}$, then there exists an $m \times m$ principal submatrix $C_S$ such that $[\epr(C_S)]_i ={\tt S}$.
\een

\end{thm}

In the rest of this paper, each instance of \
$\cdots$ is permitted to be empty.

\begin{cor}\label{NSA & ...ASN...A...}
{\rm\cite[Corollary 2.7]{EPR} ($\tt NSA$ Theorem.)}
No symmetric matrix over any field can have
$\tt{NSA}$ in its epr-sequence.
Further, no symmetric matrix over any field can have the epr-sequence $\cdots \tt{ASN} \cdots \tt{A} \cdots$.
\end{cor}

Given a matrix $B$ with a nonsingular principal
submatrix $B[\alpha]$, we denote by $B/B[\alpha]$ the
Schur complement of $B[\alpha]$ in $B$ \cite{Schur}.
The next fact is a generalization of
\cite[Proposition 2.13]{EPR} to any field;
the proof is exactly the same,
and is omitted here (we note that the proof was also omitted in \cite{EPR}).

\begin{thm}
\label{schur}
{\rm (Schur Complement Theorem.)}
Suppose $B$ is an $n \times n$
symmetric matrix over a field $\F$,
with $\rank B=r$.
Let $B[\alpha]$ be a nonsingular principal
submatrix of $B$ with $|\alpha| = k \leq r$,
and let $C = B/B[\alpha]$.
Then the following results hold.
\begin{enumerate}
\item [$(i)$]\label{p1SC} $C$ is an
$(n-k)\times (n-k)$ symmetric matrix. 
\item [$(ii)$]\label{p2SC} Assuming the indexing of $C$ is inherited from $B$, any principal minor of $C$
is given by 
\[ \det C[\gamma] = \det B[\gamma \cup \alpha]/ \det B[\alpha].\]
\item [$(iii)$]\label{p3SC} $\rank C = r-k$.
\end{enumerate}
\end{thm}

The next result, which is immediate from
the Schur Complement Theorem, has been
used implicitly in \cite{EPR} and \cite{XMR-Classif},
but we state it here in the interest of clarity
(it should be noted that this result appeared in
\cite{EPR-Hermitian} for Hermitian matrices).

\begin{cor}\label{schurAN}
{\rm (Schur Complement Corollary.)}
Let $B$ be a symmetric matrix over a field $\F$,
$\epr(B)=\ell_1 \ell_2 \cdots \ell_n$, and
let $B[\alpha]$ be a nonsingular
principal submatrix of $B$,
with $|\alpha| = k \leq \rank B$.
Let $C = B/B[\alpha]$
and $\epr(C)=\ell'_{1} \ell'_2 \cdots \ell'_{n-k}$.
Then, for $j=1, \dots, n-k$,
$\ell'_j=\ell_{j+k}$  if
$\ell_{j+k} \in \{{\tt A,N}\}$.
\end{cor}

\begin{obs}\label{Obs: Append Ns}
{\rm\cite[Observation 2.19]{EPR}}
Let $B$ be a symmetric matrix over a field $\F$, with epr-sequence $\ell_1\ell_2\cdots\ell_n$.
\ben
\item \label{Obs: Append Ns(1)}Form a matrix $B'$  from $B$ by copying the last row down and then the last column across.  Then the epr-sequence of $B'$ is $\ell_1\ell_2'\cdots\ell_n'{\tt N}$ with $\ell_i'={\tt N}$ if $\ell_i={\tt N}$ and $\ell_i'={\tt S}$ otherwise for $2\le i\le n$.
\item  \label{Obs: Append Ns(2)} Form a matrix $B''$ from $B$  by taking the direct sum with $[0]$.  Then the epr-sequence of $B''$ is $\ell_1''\ell_2''\cdots\ell_n''{\tt N}$ with $\ell_i''={\tt N}$ if $\ell_i={\tt N}$ and $\ell_i''={\tt S}$ otherwise  for $1\le i\le n$.
\een
\end{obs}

\section{Restrictions on
attainable epr-sequences over a field of characteristic $2$}\label{s: Restrictions over any field of char 2}
$\null$
\indent
Before stating our main results in
Section \ref{s: Main Results},
we devote this section towards establishing restrictions
for the attainability of epr-sequences over a field of characteristic $2$.

\begin{obs}
{\rm ({\tt NA-NS} Observation.)}
Let $B$ be a symmetric matrix over
a field of characteristic $2$, with
$\epr(B) = \ell_1 \ell_2 \cdots \ell_n$.
If
$\ell_k \ell_{k+1} = \tt NA$ or
$\ell_k \ell_{k+1} = \tt NS$
for some $k$, then
$k$ is odd and $\ell_j = \tt N$ when $j$ is odd.
\end{obs}

\bpf
Let $\pr(B) = r_0]r_1 \cdots r_n$.
Suppose
$ \ell_k \ell_{k+1} = \tt NA$ or
$ \ell_k \ell_{k+1} = \tt NS$.
Then $r_kr_{k+1} = 01$.
Since $k\geq 1$,
Theorem \ref{pr-char2 classif} implies that
$\pr(B) = 1]01\ \OL{01}\ \OL{0}$, and
therefore that $k$ is odd, and that
$\ell_j = \tt N$ when $j$ is odd.
\epf

Over a field of characteristic $2$,
the $\tt NN$ Theorem admits a generalization when
the first $\tt N$ occurs in an even position of
the epr-sequence, which is immediate from
the {\tt NA-NS} Observation and
the {\tt NN} Theorem:

\begin{obs}\label{NN generalization}
{\rm ({\tt N}-Even Observation.)}
Let $B$ be a symmetric matrix over
a field of characteristic $2$, with
$\epr(B) = \ell_1 \ell_2 \cdots \ell_n$.
Suppose $\ell_k = \tt N$ with $k$ even.
Then $\ell_j = \tt N$ for all $j \geq k$.
\end{obs}

The next observation establishes another generalization of the $\tt NN$ Theorem for epr-sequences beginning with
$\tt S$ or $\tt A$, and it is immediate from
Theorem \ref{pr-char2 classif}.

\begin{obs}\label{obs: S... & A...}
Let $B$ be a symmetric matrix over a field of characteristic $2$,
with $\epr(B) = \ell_1 \ell_2 \cdots \ell_n$.
Suppose $\ell_1 \neq \tt N$.
If $\ell_k = \tt N$ for some $k$, then
$\ell_j = \tt N$ for all $j \geq k$.
\end{obs}

In the interest of brevity,
adopting the notation in \cite{ORIGINAL},
the principal minor $\det(B[I])$ is denoted by $B_I$
(when $I = \emptyset$, $B_{\emptyset}$ is
defined to have the value 1).
Moreover, when $I = \{i_1, i_2, \dots, i_k\}$, $B_I$ is written as $B_{i_1 i_2 \cdots i_k}$.

The next result will be of
particular relevance later in this section,
and its proof resorts to Muir's law of
extensible minors \cite{Muir};
for a more recent treatment of this law,
the reader is referred to \cite{Brualdi & Schneider}.

\begin{lem}\label{Lem: X...AN implies nonzero non-principal order-(n-1) minors}
Let $n \geq 2$, and
let $B$ be a symmetric matrix over a field of characteristic $2$, with
$\epr(B) = \ell_1 \ell_2 \cdots \ell_n$.
Suppose that
$\ell_{n-1}\ell_n = \tt AN$.
Then every minor of $B$ of order $n-1$ is nonzero.
\end{lem}

\bpf
Since the desired conclusion is obvious when $n=2$, we assume that $n \geq 3$.
By hypothesis, every principal minor of $B$ of
order $n-1$ is nonzero.
Let $i,j \subseteq \{1,2, \dots, n\}$ be distinct,
and let $I = \{1,2, \dots, n\} \setminus \{i,j\}$.
Consider the $(n-1)\times (n-1)$
\emph{non-principal} submatrix resulting
from deleting row $i$ and column $j$, i.e.,
the submatrix $B[I \cup \{j\}| I \cup \{i\}]$.
Since $I$ does not contain $i$ and $j$,
using Muir's law of extensible minors
(see \cite{Muir} or \cite{Brualdi & Schneider}),
one may extend the homogenous polynomial identity
\[
B_{\emptyset}B_{ij} =
B_iB_j - \det(B[\{i\}|\{j\}])\det(B[\{j\}|\{i\}]),
\]
to obtain the identity
\[
B_{I}B_{I \cup \{i,j\}} =
B_{ I \cup \{i\} }B_{ I \cup \{j\} } -
\det(B[ I \cup \{i\}| I \cup \{j\}])
\det(B[ I \cup \{j\}|  I \cup\{i\}]).
\]
Since $B_{I \cup \{i,j\}} = \det(B)$,
and because $\ell_n = \tt N$, we must have
\[
\det(B[ I \cup \{i\}| I \cup \{j\}])
\det(B[ I \cup \{j\}|  I \cup\{i\}])=
B_{ I \cup \{i\} }B_{ I \cup \{j\} }.
\]
Then, as $\ell_{n-1} = \tt A$,
$B_{ I \cup \{i\} }B_{ I \cup \{j\} } \neq 0$,
implying that
$\det(B[ I \cup \{j\}| I \cup \{i\}]) \neq 0$.
\epf

\subsection{Restrictions on
attainable epr-sequences over $\Z_2$}\label{s: Restrictions over Z_2}
$\null$
\indent
This section focuses on establishing
restrictions for epr-sequences over $\Z_2$.

With the purpose of establishing a contrast between the attainable epr-sequences over $\Z_2$ and those over other fields of characteristic $2$, the next example exhibits matrices over a particular field of characteristic $2$ attaining epr-sequences that are not attainable over $\Z_2$
(their unattainability over $\Z_2$ is established in
this section).

\begin{ex}\label{ex: counterexamples}
{\rm
Let $\F = \Z_2$.
Consider the field
$\F[z] = \{0,1, z, z+1\}$,
where $z^2 = z+1$.
For each of the following
(symmetric) matrices over the field $\F[z]$,
$\epr(M_{\sigma}) = \sigma$,
where $\sigma$ is an epr-sequence.
{
\[
M_{\tt AAN} =
\mtx{
 1 & z & z+1 \\
 z & 1 & 0 \\
 z+1 & 0 & 1}\!\!, \
M_{\tt ASSAN} =
\mtx{
 z & 1 & z & z+1 & 0 \\
 1 & z & z+1 & 0 & 1 \\
 z & z+1 & z & 1 & z \\
 z+1 & 0 & 1 & z & z+1 \\
 0 & 1 & z & z+1 & z}\!\!, \
\]
\[
M_{\tt NANSNN} =
\mtx{
 0 & z & z+1 & 1 & 1 & 1 \\
 z & 0 & 1 & 1 & 1 & 1 \\
 z+1 & 1 & 0 & 1 & 1 & 1 \\
 1 & 1 & 1 & 0 & 1 & 1 \\
 1 & 1 & 1 & 1 & 0 & 1 \\
 1 & 1 & 1 & 1 & 1 & 0}\!\!, \
M_{\tt SAAA} =
\mtx{
 1 & 0 & 1 & 1 \\
 0 & 1 & z & z \\
 1 & z & 0 & 1 \\
 1 & z & 1 & 0}\!\!, \
\]
\[
M_{\tt SASN} =
\mtx{
 1 & z & z & 1 \\
 z & 1 & 1 & 1 \\
 z & 1 & 0 & 1 \\
 1 & 1 & 1 & 0}\!\!, \
M_{\tt SASSA} =
\mtx{
 1 & z & z & z & 1 \\
 z & 1 & 0 & 1 & 1 \\
 z & 0 & 1 & 1 & 1 \\
 z & 1 & 1 & 0 & 1 \\
 1 & 1 & 1 & 1 & 0}\!\!.
\]
}
}
\end{ex}

\begin{rem}\label{NA... and AA...}
{\rm
\indent
\ben

\item
If $B$ is an $n \times n$ symmetric matrix over $\Z_2$ having an epr-sequence starting with
$\tt AA$, then $B = I_n$.
This is because a symmetric matrix with nonzero diagonal must have each of its off-diagonal entries equal to zero in order to have all of its order-$2$ principal minors be nonzero.

\item
A similar argument shows that if an
$n \times n$ symmetric matrix $B$ over $\Z_2$ has an epr-sequence starting with $\tt NA$, then $B = A(K_n)$.
\een
}
\end{rem}

Given a sequence
$t_{i_{1}} t_{i_{2}} \cdots t_{i_{k}}$,
the notation
$\overline{t_{i_{1}} t_{i_{2}} \cdots t_{i_{k}}}$
indicates that the sequence may be repeated as many
times as desired (or it may be omitted entirely).

\begin{prop}\label{A(K_n)}
Let $n \geq 2$.
Then, over $\Z_2$,
$\epr(A(K_n)) =\tt NA\OL{NA}$ when $n$ is even, and
$\epr(A(K_n)) =\tt NA\OL{NA}N$ when $n$ is odd.
\end{prop}
\bpf
Let $\epr(A(K_n)) = \ell_1 \ell_2 \cdots \ell_n$.
Obviously, $\ell_1 = \tt N$.
Observe that, for $2 \leq q \leq n$,
every $q \times q$ principal submatrix of $B$ is
equal to $A(K_q)$.
Since $A(K_q) = J_q - I_q$,
$\det(A(K_q)) = (-1)^{q-1}(q-1)=q-1$
(in characteristic $2$).
Hence,
$\ell_q = \tt N$ when $q$ is odd and
$\ell_q = \tt A$ when $q$ is even.
\epf

\begin{lem}
{\rm ({\tt NA} Lemma.)}
Let $B$ be a symmetric matrix over $\Z_2$, with
$\epr(B) = \ell_1 \ell_2 \cdots \ell_n$.
If
$\ell_k \ell_{k+1} = \tt NA$, then
$\ell_k \cdots \ell_n = \tt NA\OL{NA}$ or
$\ell_k \cdots \ell_n = \tt NA\OL{NA}N$.
\end{lem}

\bpf
Suppose $\ell_k \ell_{k+1} = \tt NA$.
If $k=1$, then
Remark \ref{NA... and AA...} implies that
$B=A(K_n)$, and therefore that
$\epr(B) = \tt NA \OL{NA}$ or
$\epr(B) = \tt NA \OL{NA}N$
(by Proposition \ref{A(K_n)}).
Now, suppose $k\geq 2$, and that
$\ell_j \neq \tt A$ for
some even integer $j > k+1$.
By the Inheritance Theorem,
$B$ contains a singular $j \times j$
principal submatrix, $B'$, whose epr-sequence
$\ell'_1\ell'_2 \cdots \ell'_j$ has
$\ell'_k \ell'_{k+1} = \tt NA$ and $\ell'_j = \tt N$.
Since $k \geq 2$, the $\tt NN$ Theorem implies that
$\ell'_{k-1} \neq \tt N$.
Let $B'[\alpha]$ be a nonsingular
$(k-1) \times (k-1)$ principal submatrix of $B'$.
It follows from the Schur Complement Theorem that
$B'/B'[\alpha]$ is a
(symmetric) matrix of order $j-k+1$,
and from the Schur Complement Corollary that
$\epr(B'/B'[\alpha]) = \tt NA \cdots N$.
Since $\epr(B'/B'[\alpha])$ begins with $\tt NA$,
$B'/B'[\alpha] = A(K_{j-k+1})$
(by Remark \ref{NA... and AA...}).
Then, as $\epr(B'/B'[\alpha])$ ends with $\tt N$,
Proposition \ref{A(K_n)} implies that
$\epr(B'/B'[\alpha]) = \tt NA\OL{NA}N$;
hence, $j-k+1$ is odd,
which is a contradiction,
since $j$ is even and $k$ is odd.
\epf

The epr-sequence of the matrix $M_{\tt NANSNN}$
in Example \ref{ex: counterexamples} demonstrates that
the $\tt NA$ Lemma cannot be generalized to all fields of characteristic $2$.

\begin{thm}\label{AA}
{\rm ({\tt AA} Theorem.)}
If an epr-sequence containing $\tt AA$ as a
non-terminal subsequence is attainable
over $\Z_2$, then
it is the sequence $\tt \OL{A}AAA\OL{A}$.
\end{thm}

\bpf
Let $B$ be an $n \times n$ symmetric matrix over $\Z_2$, with
$\epr(B) = \ell_1 \ell_2 \cdots \ell_n$.
Suppose that $\ell_k \ell_{k+1} = \tt AA$,
where $k+1 < n$.
We now show by contradiction that
$\ell_{k+2} = \tt A$; thus, suppose
$\ell_{k+2} \neq  \tt A$.
Hence, by the Inheritance Theorem, $B$ contains a
$(k+2)\times (k+2)$ principal submatrix $C$ with
$\epr(C) = \ell'_1 \ell'_2 \cdots \ell'_{k+2}$ having
$\ell'_k \ell'_{k+1} \ell'_{k+2} = \tt AAN$.
Note that $C$ is singular.
By Remark \ref{NA... and AA...}, $k \geq 2$
(otherwise, $C = I_{3}$, which is nonsingular).
Let $I = \{1, 2, \dots, k+2\} \setminus \{1,2,3\}$.
By \cite[Theorem 2]{Hyper},
and because $C$ is over a field of characteristic $2$,
the following equation holds:
\begin{equation}\label{Hyper equation}
C_{I}^2 C_{I \cup \{1,2,3\}}^2 +
C_{I \cup \{1\}}^2 C_{I \cup \{2,3\}}^2 +
C_{I \cup \{2\}}^2 C_{I \cup \{1,3\}}^2 +
C_{I \cup \{3\}}^2 C_{I \cup \{1,2\}}^2 = 0,
\end{equation}
which is the hyperdeterminantal relation obtained
from the relation (2) appearing on \cite[p. 635]{Hyper}.
Then, as $|I| = k-1$,
the fact that
$\ell'_k \ell'_{k+1}\ell'_{k+2} = \tt AAN$
leads to a contradiction, since the quantity on the
left side of (\ref{Hyper equation}) must be nonzero.
Hence, it must be the case that $\ell_{k+2} = \tt A$.
It now follows inductively that
$\ell_k \cdots \ell_n = \tt AAA\OL{A}$.

Now,
suppose that $\ell_{j} \neq \tt A$ for some $j <k$.
Then, as $k+1 < n$, the Inverse Theorem implies that
$\epr(B^{-1})$ starts with $\tt AA$, and that
$\epr(B^{-1}) \neq \tt AA\OL{A}A$.
But, by Remark \ref{NA... and AA...}, $B^{-1} = I_n$,
implying that $\epr(B^{-1}) = \tt AA\OL{A}A$, a contradiction.

Since $\tt \OL{A}AAA\OL{A}$ is attained by $I_n$, the desired conclusion follows.
\epf

The epr-sequence of the matrix $M_{\tt AAN}$
in Example \ref{ex: counterexamples} shows that
the $\tt AA$ Theorem does not hold for all fields of characteristic $2$.

\begin{thm}\label{A...AN implies n even}
Let $n \geq 3$, and
let $B$ be a symmetric matrix over $\Z_2$, with
$\epr(B) = \ell_1 \ell_2 \cdots \ell_n$.
Suppose that $\ell_1 = \tt A$ and
$\ell_{n-1}\ell_n = \tt AN$.
Then $n$ is even.
\end{thm}

\bpf
By Lemma \ref{Lem: X...AN implies nonzero non-principal order-(n-1) minors},
every minor of $B$ of order $n-1$ is nonzero.
We claim that each row of $B$
contains an even number of nonzero entries;
to see this, let $k$ be the number of
nonzero entries of $B$ in row $i$, and
consider a calculation of $\det(B)$ via a
Laplace expansion along row $i$.
Because in the field $\Z_2$ every number is equal to
its negative, this expansion calculates $\det(B)$
by adding $k$ minors of $B$ of order $n-1$;
since each of these $k$ minors is nonzero,
and because $\det(B) = 0$,
it follows that $k$ must be even, as claimed.
Hence, the total number of nonzero entries of
$B$ must also be even.
Then, as the number of nonzero off-diagonal entries of
a symmetric matrix is always even,
it is immediate that the number of
nonzero diagonal entries of $B$ must also be even.
Finally, since the number of
nonzero diagonal entries of $B$ is $n$
(because $\ell_1 = \tt A$),
$n$ is even, as desired.
\epf

We note that the sequence
$\tt AS\OL{S}AN$ is attainable over $\Z_2$ when
its order is even (see Theorem \ref{A...}),
implying that a sequence of the form
$\tt A \cdots AN$ is not completely prohibited.
Moreover, Theorem \ref{A...AN implies n even} does not hold for all fields of characterisitic $2$ (see Example \ref{ex: counterexamples}).

In the interest of brevity when proving the
next result, define the
$n \times n$ matrix $R_{n,k}$ as follows:
For $n \geq 2$, let
\[
R_{n,k} :=
\mtx{I_{k}   &  J_{k,n-k}\\
     J_{n-k,k} &  A(K_{n-k})},
\]
where $0\leq k \leq n$
(we assume that $R_{n,k} = I_n$   when $k=n$,
and that        $R_{n,k} = A(K_{n})$ when $k=0$).

\begin{prop}\label{SA...}
An epr-sequence starting with $\tt SA$
is attainable by a symmetric matrix over $\Z_2$
if and only if
it has one of the following forms.
\[{\tt
SA\OL{SA},
\qquad
SA\OL{SA}A,
\qquad
SA\OL{SA}N.
}\]
\end{prop}

\bpf
Let $0 \leq k \leq n$ be integers.
We begin by showing that
$\det(R_{n,k}) = 0$ only when
$n$ is odd and $k$ is even.
The desired conclusion is immediate for
the case with $k=0$ (by Proposition \ref{A(K_n)}),
and, for the case with $k=n$, it is obvious
(since $R_{n,k} = I_n$).
Now, suppose $0<k< n$,
and let $C = R_{n,k}$.
Note that
$\det(C) = \det(I_{k}) \det(C/I_{k}) = \det(C/I_{k})$,
where $C/I_{k}$ is the Schur complement of
$I_{k}$ in $C$.
Then, as
\[
C/I_{k} =
A(K_{n-k}) - J_{n-k,k} \cdot J_{k,n-k} =
(1-k)J_{n-k} - I_{n-k},
\]
$\det(C) =
((1-k)(n-k) - 1)(-1)^{n-k-1} =
(k+1)n + 1$ (in characteristic $2$).
It follows that
$\det(C) = 1$   when $n$ is even, and that
$\det(C) = k$ when $n$ is odd.
We can now conclude that
$\det(R_{n,k}) = 0$ only when
$n$ is odd and $k$ is even,
as desired.

Let $\sigma$ be an epr-sequence starting with $\tt SA$.
For the first direction,
suppose that $\sigma = \epr(B)$ for some symmetric
matrix $B$ over $\Z_2$.
Let $\sigma = \ell_1 \ell_2 \cdots \ell_n$.
By hypothesis, $\ell_1\ell_2 = \tt SA$.
Without loss of generality,
suppose that the first $k$ diagonal entries of $B$
are nonzero, and suppose that the remaining
$n-k$ diagonal entries are zero.
Note that, since $\ell_1= \tt S$, $1 \leq k \leq n-1$.
It is easy to verify that the condition that
$\ell_2 = \tt A$ implies that
$B = R_{n,k}$.
It is also easy to see that for
any integer $m$ with $3 \leq m \leq n$,
any $m \times m$ principal submatrix of $R_{n,k}$
is of the form $R_{m,p}$,
where $0 \leq p \leq k$
(and $0 \leq m-p \leq n-k$).
The above argument implies that
any principal minor of $B$ of order $m$
is nonzero when $m$ is even,
implying that $\ell_j = \tt A$ whenever
$j$ is even.
Also, observe that for
any odd integer $m$ with $3 \leq m < n$,
there exists
$0 \leq p \leq k$ even,  and
$0 \leq q \leq k$ odd,
such that $R_{m,p}$ and $R_{m,q}$ are principal submatrices of $B$;
then, as
$\det(R_{m,p}) = 0$ and
$\det(R_{m,q}) \neq 0$,
$B$ contains both a
zero and a nonzero principal minor of order $m$,
implying that $\ell_j = \tt S$ whenever $j <n$ is odd.
It now follows that $B$ must have one of the desired epr-sequences.

For the other direction,
note that the order-$n$ sequence
$\tt SA\OL{SA}$ is attained by
the matrix $R_{n,1}$ when $n$ is even.
Similarly, (when $n$ is odd)
the order-$n$ sequences
$\tt SA\OL{SA}A$ and
$\tt SA\OL{SA}N$ are attained by
$R_{n,1}$ and $R_{n,2}$, respectively.
\epf

As with the previous results,
Proposition \ref{SA...} cannot be generalized either
(see Example \ref{ex: counterexamples}).

An observation following from
the {\tt NA} Lemma, the {\tt AA} Theorem and
Proposition \ref{SA...} is in order:

\begin{obs}\label{obs: XA...}
Let $B$ be a symmetric matrix over $\Z_2$, with
$\epr(B) = \ell_1 \ell_2 \cdots \ell_n$.
If $\ell_2 = \tt A$,
then $\ell_j = \tt A$ when $j$ is even.
\end{obs}

The previous and the next result also do not hold for
all fields of characteristic $2$
(see Example \ref{ex: counterexamples}).

\begin{prop}\label{SAXN}
For any $\tt X$,
the epr-sequence $\tt SAXN$
cannot occur as a subsequence of
the epr-sequence of a symmetric matrix over $\Z_2$.
\end{prop}

\bpf
Let $B$ be an $n \times n$ symmetric matrix over $\Z_2$, with
$\epr(B) = \ell_1 \ell_2 \cdots \ell_n$.
Suppose
$\ell_k \cdots \ell_{k+3} = \tt SAXN$
for some $ 1 \leq k \leq n-3$,
where $\tt X \in \{A,S,N\}$.
By Proposition \ref{SA...}, $k \geq 2$.
By the $\tt NSA$ Theorem, $\ell_{k-1} \neq \tt N$.
Let $B[\alpha]$ be a $(k-1) \times (k-1)$
nonsingular principal submatrix of $B$.
By the Schur Complement Corollary,
$\epr(B/B[\alpha]) = \tt YAZN \cdots$,
where $\tt Y, Z \in \{A, S, N\}$,
which contradicts Observation \ref{obs: XA...}.
\epf

In the epr-sequence of a symmetric matrix over a
field of characteristic \emph{not} $2$,
\cite[Theorem 2.15]{EPR} asserts that $\tt ANS$ can only
occur as the initial subsequence.
Over $\Z_2$,
the same restriction holds for $\tt ASS$:

\begin{prop}\label{ASS}
In the epr-sequence of a symmetric matrix over $\Z_2$, the subsequence $\tt ASS$ can only
occur as the initial subsequence.
\end{prop}

\bpf
Let $B$ be an $n \times n$ symmetric matrix over $\Z_2$, with
$\epr(B) = \ell_1 \ell_2 \cdots \ell_n$.
Suppose to the contrary that
$\ell_k \ell_{k+1}\ell_{k+2} = \tt ASS$
for some $ 2 \leq k \leq n-3$.
By the Inheritance Theorem,
$B$ contains a $(k+2) \times (k+2)$
principal submatrix $B'$ with
$\epr(B') = \cdots \tt XAYN$, where
$\tt X, Y \in \{A,S,N\}$.
By the {\tt NA} Lemma, $\tt X \neq N$, and,
by the {\tt AA} Theorem, $\tt X \neq A$;
hence, $\tt X = S$, so that
$\epr(B') = \cdots \tt SAYN$,
which contradicts Proposition \ref{SAXN}.
\epf

Once again, the previous result also cannot be generalized to all fields of characteristic $2$
(see Example \ref{ex: counterexamples}).

\begin{lem}\label{lem: ASA}
Let $B$ be a symmetric matrix over $\Z_2$, with
$\epr(B) = \ell_1 \ell_2 \cdots \ell_n$.
Suppose $ \ell_k \ell_{k+1} \ell_{k+2} = \tt ASA$
for some $1 \leq k \leq n-2$.
Then $\ell_1 \neq \tt N$ and the following hold.

\ben
\item If $\ell_1 = \tt A$, then $k$ is odd.
\item If $\ell_1 = \tt S$, then $k$ is even.
\een

\end{lem}

\bpf
By the $\tt NN$ Theorem and the {\tt NA-NS} Observation,
$\epr(B)$ does not begin with
$\tt NN$, $\tt NA$, nor $\tt NS$;
hence, $\ell_1 \neq \tt N$.

(1): Suppose that
$\ell_1 = \tt A$ and that $k$ is even.
Then, by the Inheritance Theorem,
$B$ contains a $(k+2) \times (k+2)$
principal submatrix $B'$
with $\epr(B') = \tt A \cdots ASA$.
By the Inverse Theorem,
$\epr(B^{-1}) = \tt SA \cdots AA$.
Since $B^{-1}$ is of order $k+2$,
Proposition \ref{SA...} implies that $k+2$ is odd,
which is a contradiction to $k$ being even.

(2): Suppose that
$\ell_1 = \tt S$ and that $k$ is odd.
Then, by the Inheritance Theorem,
$B$ contains a $(k+2) \times (k+2)$
principal submatrix $B'$
with $\epr(B') = \tt S \cdots AXA$,
where $\tt X \in \{A,S,N\}$.
Since $\tt X$ occurs in an even position,
the {\tt N}-Even Observation implies that
$\tt X \neq N$;
and, by the {\tt AA} Theorem, $\tt X \neq A$;
hence, $\tt X = S$.
By the Inverse Theorem,
$\epr((B')^{-1}) = \tt SA \cdots SA$.
Since $(B')^{-1}$ is of order $k+2$,
Proposition \ref{SA...} implies that $k+2$ is even,
a contradiction.
\epf

The inverse of the matrix $M_{\tt SASSA}$ in
Example \ref{ex: counterexamples},
whose epr-sequence is $\tt SSASA$, reveals that the previous result also cannot be generalized to
all fields of characteristic $2$;
and, for the same reasons,
the following theorem cannot be generalized either.

\begin{thm}\label{ASA anywhere}
Let $B$ be a symmetric matrix over $\Z_2$.
Suppose $\epr(B)$ contains $\tt ASA$ as a subsequence.
Then $\epr(B)$ is one of the following sequences.
\ben
\item $\tt ASA\OL{SA}$\label{ASA1};
\item $\tt ASA\OL{SA}A$\label{ASA2};
\item $\tt ASA\OL{SA}N$\label{ASA3};
\item $\tt SASA\OL{SA}$\label{ASA4};
\item $\tt SASA\OL{SA}A$\label{ASA5};
\item $\tt SASA\OL{SA}N$\label{ASA6}.
\een
\end{thm}

\bpf
Suppose that $\epr(B) = \ell_1 \ell_2 \cdots \ell_n$,
and that $\ell_k \ell_{k+1} \ell_{k+2} = \tt ASA$.
By Lemma \ref{lem: ASA}, $\ell_1 \neq \tt N$.
We proceed by examining two cases.

Case 1: $\ell_1 = \tt S$.
Because of Proposition \ref{SA...},
it suffices to show that $\ell_2 = \tt A$.
By Lemma \ref{lem: ASA}, $k$ is even.
If $k =2$, then, obviously, $\ell_2 = \tt A$.
Now, suppose $k \geq 4$.
By the Inheritance Theorem,
$B$ contains a $(k+2) \times (k+2)$
principal submatrix, $B'$, whose epr-sequence
$\ell'_1 \ell'_2 \cdots \ell'_{k+2}$ has
$\ell'_2 = \ell_2$ and
$\ell'_k \ell'_{k+1} \ell'_{k+2} =
{\tt A} \ell'_{k+1}{\tt A}$.
By the Inverse Theorem,
$\epr((B')^{-1}) =
\ell'_{k+1} {\tt A} \cdots \ell_2 \ell'_1 {\tt A}$.
It follows from Observation \ref{obs: XA...} that
$[\epr((B')^{-1})]_j = \tt A$ when $j$ is even.
Then, as $k$ is even, and because
$[\epr((B')^{-1})]_k = \ell_2$,
we must have $\ell_2 = \tt A$.

Case 2: $\ell_1 = \tt A$.
By Lemma \ref{lem: ASA}, $k$ is odd.
Let $1<j<k$ be an odd integer.
By the Inheritance Theorem,
$B$ contains a $(k+2) \times (k+2)$
principal submatrix, $B'$, whose epr-sequence
$\ell'_1 \ell'_2 \cdots \ell'_{k+2}$ has
$\ell'_j = \ell_j$ and
$\ell'_k \ell'_{k+1} \ell'_{k+2} =
{\tt A}\ell'_{k+1}{\tt A}$.
By the Inverse Theorem,
$\epr((B')^{-1}) =
\ell'_{k+1} {\tt A} \cdots \ell_j \cdots$.
It follows from Observation \ref{obs: XA...} that
$[\epr((B')^{-1})]_i = \tt A$ when $i$ is even.
Since $k+2-j$ is even,
$[\epr((B')^{-1})]_{k+2 - j} = \tt A$.
Then, as $[\epr((B')^{-1})]_{k+2 - j} = \ell'_j = \ell_j$,
we have $\ell_j = \tt A$.
We conclude that
$\ell_i = \tt A$ when
$i$ is an odd integer with  $1< i <k$.
Then, as $\ell_{k+1} = \tt S$,
the {\tt AA} Theorem implies that
$\ell_i \neq \tt A$ when
$i$ is an even integer with $1< i <k$;
and, since $\ell_k = \tt A$,
the {\tt N}-Even Observation implies that
$\ell_i \neq \tt N$ when
$i$ is an even integer with $1< i <k$.
Hence,
$\epr(B) =
{\tt ASA\OL{SA}} \ell_{k+3} \cdots \ell_n$.

If $n = k+2$, then we are done;
thus, suppose $n \geq k+3$.
Suppose to the contrary that
$\ell_q \neq \tt A$ for some odd integer $q$ with
$k+3 \leq q \leq n$.
By the Inheritance Theorem,
$B$ contains a singular $q \times q$
principal submatrix, $B'$,
whose epr-sequence
$\ell'_1 \ell'_2 \cdots \ell'_q$ has
$\ell'_i = \ell_i = \tt A$ when $i \leq k+2$ is odd,
and, obviously, $\ell'_q = \tt N$.
Let $B'[\alpha]$ be a (necessarily nonsingular)
$1 \times 1$ principal submatrix of $B'$.
By the Schur Complement Theorem,
$B'/B'[\alpha]$ is a
$(q-1) \times (q-1)$ (symmetric) matrix,
and, by the Schur Complement Corollary,
$[\epr(B'/B'[\alpha])]_2 = \tt A$, and
$[\epr(B'/B'[\alpha])]_{q-1} = \tt N$.
It follows from Observation \ref{obs: XA...}
that $[\epr(B'/B'[\alpha])]_i = \tt A$ when
$i$ is even.
Then, as $[\epr(B'/B'[\alpha])]_{q-1} = \tt N$,
$q-1$ is odd, which is a contradiction to
the fact that $q$ is odd.
We conclude that
$\ell_i = \tt A$ for all
odd $i$ with $k+3 \leq i \leq n$.

Then, as $\ell_{k+1} = \tt S$,
the {\tt AA} Theorem implies that
$\ell_i \neq \tt A$ when
$i$ is an even integer with $k+3 \leq i \leq n-1$;
and, since at least one of
$\ell_{n-1}$ and $\ell_n$ must be $\tt A$
(because one of $n-1$ and $n$ must be even)
the {\tt N}-Even Observation implies that
$\ell_i \neq \tt N$ when
$i$ is an even integer with $k+3 \leq i \leq n-1$.
It follows that
$\epr(B) = \tt ASA\OL{SA}$ when $n$ is odd,
and that either
$\epr(B) = \tt ASA\OL{SA}A$
or
$\epr(B) = \tt ASA\OL{SA}N$
when $n$ is even.
\epf

\section{Main results}\label{s: Main Results}
$\null$
\indent
In this section, a complete characterization of the epr-sequences that are attainable by a symmetric
matrix over $\Z_2$ is established.
We start by characterizing those that begin with $\tt N$.

\begin{lem}\label{NSN...NA and NSN...NAN}
Let
$M_1 = A(K_2) \oplus A(K_2) \oplus \cdots \oplus A(K_2)$ 
be $n \times n$ 
and
\[
M_2 = \mtx{M_1 & \ones_{n} \\
     \ones_{n}^T & O_1},
\]
with both matrices being over $\Z_2$.
Then
$\epr(M_1) = \tt \OL{NS}NA$ and
$\epr(M_2) = \tt \OL{NS}NAN$.
\end{lem}

\bpf
Let $\epr(M_1) = \ell_1 \ell_2 \cdots \ell_n$.
Note that $n$ is even.
The desired conclusion is obvious when $n = 2$;
hence, suppose $n \geq 4$.
It is clear that $\ell_1 \ell_2 = \tt NS$;
thus, by the {\tt NA-NS} Observation,
$\epr(M_1)$ has
$\tt N$ in every odd position.
Clearly, $M_1$ is nonsingular, implying that
$\ell_n = \tt A$.
It remains to show that $\ell_j = \tt S$ when
$j \leq n-1$ is even.
Since $\ell_n = \tt A$, by the $\tt NN$ Theorem,
$\ell_j \neq \tt N$ when $j \leq n-1$ is even.
Now, because of the {\tt NA} Lemma,
to show that $\ell_j \neq \tt A$ when
$j \leq n-1$ is even,
it suffices to show that $\ell_{n-2} = \tt S$.
Clearly, $M_1(\{2,4\})$ is singular
(since it contains a zero row).
Then, as $\ell_{n-2} \neq \tt N$
(because $n-2$ is even),
$\ell_{n-2} = \tt S$.

Let $\epr(M_2) = \ell'_1 \ell'_2 \cdots \ell'_{n+1}$.
The assertion is clear when $n = 2$
(note that $n$ is even, and that
$M_2$ is of order $n+1$, not $n$);
hence, suppose $n \geq 4$.
Since (clearly) $\ell'_1 \ell'_2 = \tt NS$,
the {\tt NA-NS} Observation implies that $\epr(M_2)$ has
$\tt N$ in every odd position.
Since $M_1$ is a principal submatrix of $M_2$,
and because $\epr(M_1) = \tt NS\OL{NS}NA$,
it is immediate that
$\epr(M_2) = {\tt NS\OL{NS}N}\ell'_n \tt N$.
We now show that $\ell'_n = \tt A$.
Observe that any $n \times n$ principal submatrix
of $M_2$ is either $M_1$, which is nonsingular, or is one that is permutationally similar  to the matrix
\[
C = \mtx{C(\{n\}) & \ones_{n-1} \\
     \ones_{n-1}^T & O_1},
\]
where
$C(\{n\}) = O_1 \oplus A(K_2) \oplus A(K_2) \oplus \cdots
\oplus A(K_2)$.
Let $C'$ be the matrix obtained from $C$ by first
subtracting its first row from
rows $2, 3, \dots, n-1$,
and then subtracting the first column of
the resulting matrix from
columns $2, 3, \dots, n-1$.
Now observe that
$\det(C') = -\det(C'(\{1,n\}))$, where
$C'(\{1,n\}) =
A(K_2) \oplus A(K_2) \oplus \cdots \oplus A(K_2)$,
which is a nonsingular matrix (of order $(n-2)$).
Hence, $\det(C') \neq 0$.
Then, as $\det(C) = \det(C')$,
$C$ is nonsingular.
We conclude that $\ell'_n = \tt A$.
\epf

\begin{thm}\label{N...}
An epr-sequence starting with $\tt N$ is attainable
by a symmetric matrix over $\Z_2$
if and only if
it has one of the following forms:
\ben
\item $\tt NA\OL{NA}$;
\item $\tt NA\OL{NA}N$;
\item $\tt \OL{NS}N\OL{N}$;
\item $\tt NS\OL{NS}NA$;
\item $\tt NS\OL{NS}NAN$.
\een
\end{thm}

\bpf
Let $\sigma = \ell_1 \ell_2 \cdots \ell_n$ be
an epr-sequence with $\ell_1 = \tt N$.
Suppose that $\sigma = \epr(B)$, where
$B$ is a symmetric matrix over $\Z_2$.
If $n=1$, then
$\sigma = \tt \OL{NS}N\OL{N}$ with
$\tt \OL{NS}$ and $\tt \OL{N}$ empty.
Suppose $n \geq 2$.
If $\ell_2 = \tt N$, then,
by the $\tt NN$ Theorem,
$\sigma = \tt \OL{NS}NN\OL{N}$
with $\tt \OL{NS}$ empty.
If $\ell_2 = \tt A$, then,
by the {\tt NA} Lemma,
$\sigma = \tt NA\OL{NA}$ or
$\sigma = \tt NA\OL{NA}N$.

Finally, suppose $\ell_2 = \tt S$.
Since an attainable epr-sequence cannot
end in $\tt S$, $n \geq 3$.
By the {\tt NA-NS} Observation,
$\ell_j = \tt N$ when $j$ is odd.
Hence, $\rank(B)$ is even.
We now show that $\tt SNA$ cannot
occur as a subsequence of
$\ell_1 \ell_2 \cdots \ell_{n-2}$.
Suppose to the contrary that
$\ell_{k-1} \ell_{k} \ell_{k+1} = \tt SNA$,
where $3 \leq k \leq n-3$.
Clearly, since $\ell_j = \tt N$ when $j$ is odd,
$k$ is odd and $\ell_{k+2} = \tt N$.
By the Inheritance Theorem,
$B$ contains a $(k+3) \times (k+3)$
principal submatrix $B'$ with
$\epr(B') = \cdots \tt SNANX$, where
$\tt X \in \{A, N\}$.
If $\tt X = A$, then, by the Inverse Theorem,
$\epr((B')^{-1}) = \tt NANS \cdots$,
which contradicts the {\tt NA} Lemma.
Hence, $\tt X = N$, and therefore
$\epr(B') = \cdots \tt SNANN$,
which contradicts the {\tt NA} Lemma.
We conclude that $\tt SNA$ cannot
occur as a subsequence of
$\ell_1 \ell_2 \cdots \ell_{n-2}$.
Now, let $r = \rank(B)$;
hence, $\ell_r \neq \tt N$.
Then, as $r$ is even, $\ell_{r-1} = \tt N$
(because $r-1$ is odd).
Since $\ell_j = \tt N$ when $j$ is odd,
the $\tt NN$ Theorem implies that
$\ell_i \neq \tt N$ when $i \leq r-1$ is even.
We proceed by considering two cases.

Case 1: $r \geq n-1$.
First, suppose $r = n-1$.
Since $r$ is even, $r+1 = n$ is odd,
implying that $\ell_n = \tt N$.
Hence,
$\ell_{n-1}\ell_n = \tt AN$ or
$\ell_{n-1}\ell_n = \tt SN$.
Then, as $\ell_2 = \tt S$,
and because $\tt SNA$ cannot occur as a
subsequence of $\ell_1 \ell_2 \cdots \ell_{n-2}$,
it follows inductively that
$\sigma = \tt NS\OL{NS}NAN$ or
$\sigma = \tt \OL{NS}NSN$.
Now, suppose $r=n$;
hence, $n$ is even and $\ell_n = \tt A$.
Since $\ell_{r-1} = \tt N$,
$\ell_{n-1}\ell_n = \tt NA$.
Then, as $\ell_2 = \tt S$,
and because $\tt SNA$ cannot occur as a
subsequence of $\ell_1 \ell_2 \cdots \ell_{n-2}$,
it follows inductively that
$\sigma = \tt NS\OL{NS}NA$.

Case 2: $r \leq n-2$.
Hence,
$\ell_{r+1} \cdots \ell_n = \tt NN\OL{N}$.
Since $\ell_{r-1} = \tt N$ and $\ell_r \neq \tt N$,
$\ell_{r-1} \cdots \ell_n = \tt NANN\OL{N}$ or
$\ell_{r-1} \cdots \ell_n = \tt NSNN\OL{N}$;
but the former case contradicts the {\tt NA} Lemma,
implying that $\ell_{r-1} \cdots \ell_n = \tt NSNN\OL{N}$.
Then, as $\ell_2 = \tt S$,
and because $\tt SNA$ cannot occur as a
subsequence of $\ell_1 \ell_2 \cdots \ell_{n-2}$,
it follows inductively that
$\sigma = \tt NS\OL{NS}NN\OL{N}$.


For the other direction, we show that
each of the sequences listed above is attainable.
Assume that the sequence under
consideration has order $n$.
The sequences
$\tt NA\OL{NA}$ and $\tt NA\OL{NA}N$ are
attainable by Proposition \ref{A(K_n)}.
When $\tt \OL{NS}$ is non-empty
the sequence $\tt \OL{NS}N\OL{N}$ is attainable
by applying
Observation \ref{Obs: Append Ns}(2)
to the sequence $\tt NA\OL{NA}$;
and, when $\tt \OL{NS}$ is empty, it is attained by $O_n$.
Finally, the sequences
$\tt NS\OL{NS}NA$ and $\tt NS\OL{NS}NAN$
are attainable by Lemma \ref{NSN...NA and NSN...NAN}.
\epf

Naturally, due to the dependence of Theorem \ref{N...} on the results of Section \ref{s: Restrictions over Z_2},
this theorem does not hold for other fields.

Some lemmas are necessary before stating
the second of our three main results in
Theorem \ref{A...}.

\begin{lem}\label{ASS...A and ASS...AA}
Let $n \geq 4$, $m \geq 5$, and let
\[
M_{\tt ASA} =
\mtx{I_2& \ones_2 \\
     \ones_2^T & J_1},\
M_{\tt ASAA} =
\mtx{I_2 & J_2 \\
     J_2 & I_2}
\]
be over $\Z_2$.
Let
$B = I_{n-3} \oplus M_{\tt ASA}$,
$B' = I_{m-4} \oplus M_{\tt ASAA}$,
$\epr(B) = \ell_1 \ell_2 \cdots \ell_n$ and
$\epr(B') = \ell'_1 \ell'_2 \cdots \ell'_m$.
Then
$\epr(M_{\tt ASA}) = \tt ASA$,
$\epr(M_{\tt ASAA}) = \tt ASAA$,
$\ell_1\ell_2\ell_3 = \ell'_1\ell'_2\ell'_3 = \tt ASS$,
$\ell_{n-1}\ell_n = \tt SA$ and
$\ell'_{m-1} \ell'_m = \tt AA$.
\end{lem}

\bpf
All of the assertions above are easily verified,
except $\ell'_{m-1} = \tt A$, which we now prove.
The case with $m=5$ is easy to check;
thus, suppose $ m \geq 6$.
Note that,
since every $3 \times 3$ principal submatrix of
the ($4 \times 4$) matrix $M_{\tt ASAA}$ is nonsingular,
and because every $(m-5) \times (m-5)$ principal
submatrix of $I_{m-4}$ is also nonsingular,
deleting row $i$ and column $i$ of $B'$
results in a matrix
that is a direct sum of two nonsingular matrices;
hence, every $(m-1) \times (m-1)$
principal submatrix of $B'$ is nonsingular,
implying that $\ell'_{m-1} = \tt A$.
\epf

A matrix that will play an important
role here is defined as follows:
For $n \geq 2$, let $F_{n}$ be the
$n \times n$ matrix resulting from
replacing the first diagonal entry of $A(K_n)$ with $1$.

\begin{lem}\label{F_n is nonsigular}
Let $n \geq 2$, and let $F_n$ be over $\Z_2$.
Then $F_n$ is nonsingular.
\end{lem}

\bpf
The assertion is obvious when $n=2$;
thus, assume $n \geq 3$.
Observe that
\[
\det(F_n) =
\det(F_n[\{1\}])\det(F_n/F_n[\{1\}])=
\det(J_1)\det(F_n/J_1)=
\det(F_n/J_1),
\]
where
\[
F_n/J_1 =
F_n[\{2,\dots,n\}] -
\ones_{n-1} \cdot (J_1)^{-1} \cdot \ones_{n-1}^T =
A(K_{n-1}) - J_{n-1} = -I_{n-1}.
\]
Hence, $\det(F_n) = \det(-I_{n-1}) \neq 0$.
\epf

\begin{lem}\label{ASS...SAN, n=4k+2}
Let $n = 4k+2$, where $k \geq 1$ is an integer.
Let $m= \frac{n}{2}$,
let
\[
B = \mtx{J_{m} &   I_{m} \\
         I_{m} &   I_{m}}
\]
be over $\Z_2$,
and let $\epr(B) = \ell_1 \ell_2 \cdots \ell_n$.
Then
$\ell_1 \ell_2 \ell_3 = \tt ASS$ and
$\ell_{n-1}\ell_n = \tt AN$.
\end{lem}

\bpf
It is easily seen that
$\ell_1 \ell_2 \ell_3 = \tt ASS$.
Next, we show that $\ell_n = \tt N$.
Observe that
$\det(B) =
\det(I_{m})\det(B/I_{m})$,
where
\[
B/I_{m} =
J_{m} -
I_{m}\cdot
(I_{m})^{-1}
\cdot I_{m} = A(K_{m}).
\]
Since $m=\frac{n}{2} = 2k+1$ is odd,
Proposition \ref{A(K_n)} implies that
$A(K_{m})$ is singular, implying that
$B/I_{m}$ is singular, and therefore that
$\det(B) = 0$;
hence, $\ell_n = \tt N$.

Now, to see that $\ell_{n-1} = \tt A$, note that the
$(n-1)\times (n-1)$ principal submatrix resulting from the  deletion of the $i$th row and $i$th column of $B$ must be one of the following two matrices:
\[
C_1 = \mtx{J_{m-1} & X^T    \\
           X     & I_{m}}, \
C_2 = \mtx{J_{m}   & X       \\
           X^T   & I_{m-1}},
\]
where
$X=I_m(\emptyset, \{q\})$
and $q \in \{1,2, \dots, m\}$
is the unique integer such that
$i=q$ or $i= m+q$
(that is,
$q=i$ if $1 \leq i \leq m$, and
$q=i-m$ if $m+1 \leq i \leq n$).
Observe that
$\det(C_1) = \det(I_m)\det(C_1/I_m)$ and
$\det(C_2) = \det(I_{m-1})\det(C_2/I_{m-1})$,
where
$C_1/I_m = J_{m-1}-X^TX$ and
$C_2/I_{m-1} = J_{m}-XX^T$
are the Schur complements of
$I_m$ and $I_{m-1}$ in
$C_1$ and $C_2$, respectively.
Since $X^TX = I_{m-1}$, $C_1/I_m = A(K_{m-1})$.
Then, as $m-1 = 2k$ is even,
Proposition \ref{A(K_n)} implies that
$C_1/I_m$ is nonsingular;
hence, $\det(C_1) \neq 0$.

Finally, observe that
$XX^T$ is the
$m \times m$ matrix resulting from replacing the
$q$th diagonal entry of $I_m$ with $0$.
Hence, $C_2/I_{m-1}$ is the matrix resulting
from replacing the $q$th diagonal entry of
$A(K_m)$ with $1$.
Then, as $C_2/I_{m-1}$ is permutationally similar to
the nonsingular matrix $F_m$
(see Lemma \ref{F_n is nonsigular}),
$C_2/I_{m-1}$ is nonsingular,
implying that $\det(C_2) \neq 0$.
\epf

A worthwhile observation is that
the condition that $n$ is equal to $2$ modulo $4$ in
Lemma \ref{ASS...SAN, n=4k+2}
was of relevance when showing that $\det(B) = 0$,
as it is consistent with the proof of
Theorem \ref{A...AN implies n even}, from which
it can be deduced that, in order to have $\det(B) = 0$,
it is necessary for $B$ to
contain an even number of nonzero entries in each row
(observe that $B$ contains $\frac{n}{2} +1=2(k+1)$ nonzero entries in each of the first $\frac{n}{2}$ rows,
and 2 nonzero entries in each of the remaining rows).
For the same reasons, the congruence modulo $4$ of $n$ in
the following lemma will once again be of relevance.

\begin{lem}\label{ASS...SAN, n=4k}
Let $n = 4k$, where $k \geq 2$ is an integer.
Let $m = \frac{n}{2}$,
let
\[
B =
\mtx{J_{m-1} & W     \\
     W^T     & I_{m+1} }
\]
be over $\Z_2$,
where $W=[I_{m-1}, J_{m-1,2}]$,
and let $\epr(B) = \ell_1 \ell_2 \cdots \ell_n$.
Then
$\ell_1 \ell_2 \ell_3 = \tt ASS$ and
$\ell_{n-1}\ell_n = \tt AN$.
\end{lem}

\bpf
It is easily verified that
$\ell_1 \ell_2 \ell_3 = \tt ASS$.
Now we verify that $\ell_n = \tt N$.
Observe that
$\det(B) = \det(I_{m+1})\det(B/I_{m+1})$, where
$B/I_{m+1}$ is the Schur complement of
$B[\{m, m+1, \dots, n\}] = I_{m+1}$ in $B$.
Note that
\[
B/I_{m+1} =
J_{m-1} - WW^T =
J_{m-1} - I_{m-1} - 2J_{m-1}=
A(K_{m-1}) - 2J_{m-1}.
\]
Hence, $B/I_{m+1} = A(K_{m-1})$ (in characteristic $2$).
Then, as $m-1 = 2k-1$ is odd,
Proposition \ref{A(K_n)} implies that
$B/I_{m+1}$ is singular.
It follows that $\det(B) = 0$, and
therefore that $\ell_n = \tt N$.

Now we show that $\ell_{n-1} = \tt A$.
Let
$\alpha_1 = \{1, 2, \dots, m-1\}$,
$\alpha_2 = \{m, m+1, \dots, n-2\}$ and
$\alpha_3 = \{n-1, n\}$.

Let $B'$ be the matrix obtained from $B$ by
deleting its $i$th row and $i$th column.
Let $q=i-(m-1)$.
Suppose that
$M_1 = B'$ if $i \in \alpha_1$, that
$M_2 = B'$ if $i \in \alpha_2$, and that
$M_3 = B'$ if $i \in \alpha_3$.
It is easy to see that
\[
C_1 = \mtx{J_{m-2} & X \\
            X^T     & I_{m+1}},\
C_2 = \mtx{J_{m-1}  & Y \\
            Y^T      & I_{m}},\
C_3 = \mtx{J_{m-1} & Z \\
            Z^T     & I_{m}},
\]
where
\[
X = [I_{m-1}(\{i\},\emptyset), J_{m-2,2}],\
Y = [I_{m-1}(\emptyset, \{q\}), J_{m-1,2}],\
Z = [I_{m-1}, \ones_{m-1}].
\]
We proceed to show that $B'$ is nonsingular
by considering the three cases outlined above.

Case 1: $B' = C_1$.
Note that
$\det(C_1) = \det(I_{m+1}) \det(C_1/I_{m+1})$,
where $C_1/I_{m+1}$ is the Schur Complement of
$I_{m+1}$ in $C_1$,
and that
\[
C_1/I_{m+1} =
J_{m-2} - XX^T =
J_{m-2} - I_{m-2}-2J_{m-2} =
A(K_{m-2})-2J_{m-2}.
\]
Hence, $C_1/I_{m+1} = A(K_{m-2})$ (in characteristic $2$).
Then, as $m-2 = 2k-2$ is even,
Proposition \ref{A(K_n)} implies that
$C_1/I_{m+1}$ is nonsingular, implying that
$\det(C_1) \neq 0$.

Case 2: $B' = C_2$.
Then
$\det(C_2) = \det(I_{m}) \det(C_2/I_{m})$,
where $C_2/I_{m}$ is the Schur complement of
$I_{m}$ in $C_2$, and
\[
C_2/I_{m} =
J_{m-1} - YY^T =
J_{m-1} -
I_{m-1}(\emptyset, \{q\}) \cdot I_{m-1}(\{q\}, \emptyset) - 2J_{m-1}.
\]
Note that
$I_{m-1}(\emptyset,\{q\})\cdot I_{m-1}(\{q\},\emptyset)$
is the matrix obtained from $I_{m-1}$ by replacing its $q$th diagonal entry with $0$.
Then, as $2J_{m-1}=O_{m-1}$ (in characteristic $2$),
$C_2/I_{m}$ is the matrix obtained from $A(K_{m-1})$ by
replacing its $q$th diagonal entry with $1$.
Hence, $C_2/I_{m}$ is permutationally similar to
the nonsingular matrix $F_{m-1}$
(see Lemma \ref{F_n is nonsigular}).
It follows that $C_2/I_{m}$ is nonsingular,
and therefore that $\det(C_2) \neq 0$.

Case 3: $B' = C_3$.
Then
$\det(C_3) = \det(I_{m}) \det(C_3/I_{m})$,
where $C_3/I_{m}$ is the Schur complement of
$I_{m}$ in $C_3$, and
\[
C_3/I_{m} =
J_{m-1} - ZZ^T =
J_{m-1} - I_{m-1} - J_{m-1}
= -I_{m-1}.
\]
It follows that $C_3/I_{m}$ is nonsingular,
and therefore that $\det(C_3) \neq 0$.
\epf

\begin{lem}\label{ASA... matrices}
The following epr-sequences are attainable by a
symmetric matrix over $\Z_2$.
\[{\tt
ASA\OL{SA},
\qquad
ASA\OL{SA}A,
\qquad
ASA\OL{SA}N.
}\]
\end{lem}

\bpf
The attainability of $\tt ASA\OL{SA}$ follows
by observing that, by the Inverse Theorem,
the inverse of any (symmetric) matrix attaining the sequence
$\tt SA\OL{SA}A$, which is attainable by
Proposition \ref{SA...}, has epr-sequence
$\tt ASA\OL{SA}$.

Now, for $n \geq 4$ even, we show that the matrix
\[
B=\mtx{I_{n-2} & J_{n-2, 2}\\
     J_{2,n-2} & I_2}
\]
has epr-sequence $\tt ASA\OL{SA}A$.
Because of Theorem \ref{ASA anywhere},
it suffices to show that
$\epr(B)$ begins with $\tt ASA$,
and that it ends with $\tt A$.
Observe that
$\det(B) = \det(I_{n-2}) \det(B/I_{n-2})$,
where
\[B/I_{n-2} =
I_2 - J_{2,n-2} \cdot (I_{n-2})^{-1} \cdot J_{n-2,2} =
I_2 - (n-2)J_2.
\]
Since $n$ is even,
$B/I_{n-2} = I_2$ (in characteristic $2$);
hence, $B$ is nonsingular.
It is clear that $\epr(B)$ begins with $\tt AS$.
Finally, note that each $3 \times 3$
principal submatrix of $B$ must be $I_3$ or
one of the following:
\[
\mtx{I_2& \ones_2 \\
     \ones_2^T & J_1},
\qquad
\mtx{J_1& \ones_2^T \\
     \ones_2 & I_2}.
\]
Then, as each of these
$3 \times 3$ matrices is nonsingular,
$\epr(B)$ begins with $\tt ASA$, as desired.


With $n \geq 5$, let $B$ be an $n \times n$ symmetric matrix with epr-sequence $\tt SASA\OL{SA}N$,
which is attainable by Proposition \ref{SA...}.
Note that $n$ is odd.
Let $\alpha \subseteq \{1,2, \dots, n\}$ with
$|\alpha| = 1$ be such that $B[\alpha]$ is nonsingular.
Let $\epr(B/B[\alpha]) =
\ell'_1 \ell'_2 \cdots \ell'_{n-1}$.
We now show that
$\epr(B/B[\alpha]) = \tt ASA\OL{SA}N$.
By the Schur Complement Corollary,
$\ell'_j = \tt A$ when $j$ is odd, and
$\ell'_{n-1} = \tt N$.
Since $n-2$ is odd, $\ell'_{n-2} = \tt A$.
Since $\ell'_{n-1} = \tt N$,
the {\tt AA} Theorem implies that
$\ell'_j \neq \tt A$ when $j \leq n-3$ is even.
Finally, as $\ell'_{n-2} = \tt A$,
the {\tt N}-Even Observation implies that
$\ell'_j = \tt S$ when $j \leq n-3$ is even.
It follows that
$\epr(B/B[\alpha]) = \tt ASA\OL{SA}N$, as desired.
\epf

Before stating our characterization of
the epr-sequences that begin with
$\tt A$ in the next theorem,
something needs to be clarified:
\cite[Corollary 2.22]{EPR} claims that the sequence
$\tt AS\OL{S}AAA\OL{A}$ is attainable over $\Z_2$;
this claim is false:
Observe that it contradicts the {\tt AA} Theorem.
But it should be noted that
\cite[Corollary 2.22]{EPR} becomes true once
the field is restricted to be of characteristic $0$,
since it relies on \cite[Proposition 2.18]{EPR}.

\begin{thm}\label{A...}
An epr-sequence of order $n$,
and starting with $\tt A$, is attainable
by a symmetric matrix over $\Z_2$
if and only if
it has one of the following forms:
\ben
\item $\tt A\OL{A}$\label{A1};
\item $\tt A\OL{S}N\OL{N}$\label{A2};
\item $\tt ASS\OL{S}A$\label{A3};
\item $\tt ASS\OL{S}AA$\label{A4};
\item $\tt ASSS\OL{S}AN$ with $n$ even\label{A5};
\item $\tt ASA\OL{SA}$\label{A6};
\item $\tt ASA\OL{SA}A$\label{A7};
\item $\tt ASA\OL{SA}N$\label{A8}.
\een
\end{thm}

\bpf
Let $\sigma = \ell_1 \ell_2 \cdots \ell_n$ be
an epr-sequence with $\ell_1 = \tt A$.
Suppose that $\sigma = \epr(B)$, where
$B$ is a symmetric matrix over $\Z_2$.
If $n=1$ or $n=2$, then $\sigma$ is
$\tt A$, $\tt AA$, or $\tt AN$,
all of which are listed above.
Suppose $n \geq 3$.
If $\ell_2 = \tt A$ or $\ell_2 = \tt N$, then
the {\tt AA} Theorem and the
{\tt N}-Even Observation imply that
$\sigma$ is either $\tt AAA\OL{A}$ or $\tt ANN\OL{N}$.
Now, suppose $\ell_2 = \tt S$.
If $\sigma$ contains the subsequence $\tt ASA$, then,
by Theorem \ref{ASA anywhere},
$\sigma$ is either
$\tt ASA\OL{SA}$,
$\tt ASA\OL{SA}A$, or
$\tt ASA\OL{SA}N$.
Now, suppose $\sigma$ does not contain $\tt ASA$.
Hence,
$\ell_3 = \tt N$ or
$\ell_3 = \tt S$, and
$n \geq 4$.
If $\ell_3 = \tt N$, then
Observation \ref{obs: S... & A...} implies that
$\sigma = \tt ASNN\OL{N}$.
Now, assume that $\ell_3 = \tt S$.
Let $k$ be a minimal integer with $3 \leq k \leq n-1$
such that
$\ell_k \ell_{k+1} = \tt SN$ or
$\ell_k \ell_{k+1} = \tt SA$.
Hence,
$\ell_1\ell_2 \cdots \ell_k = \tt ASS\OL{S}$.
If $\ell_{k+1} = \tt N$, then
Observation \ref{obs: S... & A...}
implies that $\sigma = \tt ASS\OL{S}N\OL{N}$.
Now, assume that
$\ell_{k+1} = \tt A$.
If $n = k+1$, then
$\sigma = \tt ASS\OL{S}A$.
Thus, suppose $n \geq k+2$.

We now show that $n = k+2$.
Suppose to the contrary that $n \geq k+3$.
By the {\tt AA} Theorem,
$\ell_{k+2} \neq \tt A$.
If $\ell_{k+2} = \tt N$, then
Observation \ref{obs: S... & A...} implies that
$\sigma$ contains $\tt SANN$,
which is prohibited by Proposition \ref{SAXN};
hence, $\ell_{k+2} = \tt S$, so that
$\ell_k \ell_{k+1}\ell_{k+2} = \tt SAS$.
Then, as $\sigma$ does not contain $\tt ASA$,
and because $\tt SASN$ is prohibited by
Proposition \ref{SAXN},
$\ell_{k+3} = \tt S$, implying that
$\sigma$ contains $\tt ASS$ as a
non-initial subsequence,
which contradicts Proposition \ref{ASS}.
It follows that $n = k+2$, and
therefore that $\sigma$ is either
$\tt ASS\OL{S}AA$ or $\tt ASS\OL{S}AN$;
in the case with $\sigma = \tt ASS\OL{S}AN$,
Theorem \ref{A...AN implies n even} implies that
$n$ is even, and therefore that
$\sigma = \tt ASSS\OL{S}AN$.

Now, we establish the other direction.
As before, we assume that
the sequence under consideration has order $n$.
First, the sequence $\tt A\OL{A}$ is attained by $I_n$.
The sequence $\tt A\OL{S}N\OL{N}$ is attainable
by applying Observation \ref{Obs: Append Ns}(1) to the
sequence $\tt A \OL{A}$.
To see that
$\tt ASS\OL{S}A$ and $\tt ASS\OL{S}AA$ are attainable,
observe that the matrices $B$ and $B'$ in
Lemma \ref{ASS...A and ASS...AA}
must attain these sequences, respectively,
since the epr-sequence
of these matrices must be one of those listed above.
Similarly, when $n$ is even,
one of the two matrices in the statements of
Lemma \ref{ASS...SAN, n=4k+2} and
Lemma \ref{ASS...SAN, n=4k}
is required to attain the sequence
$\tt ASSS\OL{S}AN$.
Finally, the attainability of
$\tt ASA\OL{SA}$,
$\tt ASA\OL{SA}A$, and
$\tt ASA\OL{SA}N$ follows from
Lemma \ref{ASA... matrices}.
\epf

The reader is once again referred to
Example \ref{ex: counterexamples} to see why Theorem \ref{A...} cannot be generalized to other fields.

As before, we need more lemmas in order
to prove the last of our three main results.

For an integer $n \geq 2$ and $k \in \{1,2, \dots, n\}$,
we let $e_k^n$ denote the column vector of length $n$ with
the $k$th entry equal to $1$ and
every other entry equal to zero;
moreover, let
\[
G_n := \mtx{J_1      &  (e_1^{n-1})^T \\
            e_1^{n-1}  &  F_{n-1}      }.
\]

\begin{lem}\label{SS...AN, n even}
Let $n \geq 4$ be an even integer, let
$G_n$ be over $\Z_2$,
and let $\epr(G_n) = \ell_1 \ell_2 \cdots \ell_n$.
Then
$\ell_1 \ell_2 = \tt SS$ and
$\ell_{n-1} \ell_{n} = \tt AN$.
\end{lem}

\bpf
It is easily verified that $\ell_1 \ell_2 = \tt SS$.
The final assertion is easy to check when $n = 4$;
thus, suppose $n \geq 5$.
Observe that any $(n-1) \times (n-1)$ principal submatrix of $G_n$ has one of the following forms:
$G_{n-1}$, $F_{n-1}$ or $J_1 \oplus A(K_{n-2})$.
Hence,
to show that $\ell_{n-1} \ell_{n} = \tt AN$,
it suffices to show that
$G_{n-1}$, $F_{n-1}$ and $J_1 \oplus A(K_{n-2})$
are nonsingular, and that $G_n$ is singular.
By Lemma \ref{F_n is nonsigular},
$F_{n-1}$ is nonsingular.
By Proposition \ref{A(K_n)},
and because $n-2$ is even,
$J_1 \oplus A(K_{n-2})$ is also nonsingular.

Finally, we show that
$\det(G_{n-1}) \neq 0$ and $\det(G_{n}) = 0$.
Let $q \in \{n-1, n\}$.
Observe that
$\det(G_q) = \det(F_{q-1}) - \det(A(K_{q-2}))$.
Since $\det(F_{q-1}) \neq 0$,
$\det(G_q)= 1 - \det(A(K_{q-2}))$ (in characteristic $2$).
Hence, $\det(G_q) = 0$ if and only if
$\det(A(K_{q-2})) \neq 0$.
It follows from Proposition \ref{A(K_n)} that
$\det(G_q) = 0$ if and only if $q$ is even.
Then, as $n$ is even,
$\det(G_{n-1}) \neq 0$ and $\det(G_{n}) = 0$.
\epf

\begin{lem}\label{SS...AN, n odd}
Let $n \geq 5$ be an odd integer.
Then there exists a symmetric matrix over $\Z_2$ whose epr-sequence
$\ell_1 \ell_2 \cdots \ell_n$ has
$\ell_1 \ell_2 = \tt SS$ and
$\ell_{n-1}\ell_n = \tt AN$.
\end{lem}

\bpf
Clearly, $n+1$ is even and $n+1 \geq 6$.
Let $m=\frac{n+1}{2}$, and let
\[
B'=\mtx{J_m & I_m \\
       I_m & I_m}, \
B''=\mtx{J_{m-1} & W \\
        W^T     & I_{m+1}},
\]
where $W=[I_{m-1},J_{m-1,2}]$.
Observe that $B'$ and $B''$ are
$(n+1)\times(n+1)$ symmetric matrices.
Let
$\epr(B') = \ell'_1\ell'_2 \cdots \ell'_{n+1}$ and
$\epr(B'') = \ell''_1\ell''_2 \cdots \ell''_{n+1}$.
We consider two cases:

Case 1: $n+1=4k+2$ for some integer $k \geq 1$.
Observe that, by Lemma \ref{ASS...SAN, n=4k+2},
$\ell'_{n}\ell'_{n+1} = \tt AN$.
Let $\alpha = \{n+1\}$,
let $C = B'/B'[\alpha]$, and
let $\epr(C) = \ell_1\ell_2 \cdots \ell_n$.
We now show that $C$ is a matrix with
the desired properties.
By the Schur Complement Corollary,
$\ell_{n-1}\ell_n = \tt AN$.
To show that $\ell_1\ell_2= \tt SS$,
first, observe that,
by the Schur Complement Theorem,
and because $\det(B'[\alpha]) = 1$
(in characteristic $2$),
\[
\det(C[\{n\}]) = \det(B'[\{n\} \cup \alpha]),\ \
\det(C[\{m\}]) = \det(B'[\{m\} \cup \alpha]),
\]
\[
\det(C[\{n-1,n\}]) = \det(B'[\{n-1,n\} \cup \alpha]),\ \
\det(C[\{1,2\}]) = \det(B'[\{1,2\} \cup \alpha]).
\]
Then, by observing that
$B'[\{n\} \cup \alpha]=I_2$, that
$B'[\{m\} \cup \alpha]= J_2$, that
$B'[\{n-1,n\} \cup \alpha]= I_3$ and that
$B'[\{1,2\} \cup \alpha]= J_2 \oplus J_1$,
we conclude that
$\det(C[\{n\}])$ and $\det(C[\{n-1,n\}])$ are nonzero,
and that
$\det(C[\{m\}])$ and $\det(C[\{1,2\}])$
are zero.
Hence, $\ell_1 \ell_2 = \tt SS$.

Case 2: $n+1  = 4k$ for some integer $k \geq 2$.
Observe that, by Lemma \ref{ASS...SAN, n=4k},
$\ell''_{n}\ell''_{n+1} = \tt AN$.
Let $\alpha = \{n+1\}$,
let $C = B''/B''[\alpha]$, and
let $\epr(C) = \ell_1\ell_2 \cdots \ell_n$.
As in Case 1, we show that $C$ is a matrix with
the desired properties.
By the Schur Complement Corollary,
$\ell_{n-1}\ell_n = \tt AN$.
To show that $\ell_1\ell_2= \tt SS$,
first, observe that,
by the Schur Complement Theorem,
and because $\det(B''[\alpha]) = 1$
(in characteristic $2$),
\[
\det(C[\{n\}]) = \det(B''[\{n\} \cup \alpha]),\ \
\det(C[\{1\}]) =\det(B''[\{1\} \cup \alpha]),
\]
\[
\det(C[\{n-1,n\}]) = \det(B''[\{n-1,n\} \cup \alpha]),\ \
\det(C[\{1,2\}]) = \det(B''[\{1,2\} \cup \alpha]).
\]
Then, by observing that
$B''[\{n\} \cup \alpha]=I_2$,
$B''[\{1\} \cup \alpha]= J_2$,
$B''[\{n-1,n\} \cup \alpha]= I_3$ and
$B''[\{1,2\} \cup \alpha]= J_3$,
we conclude that
$\det(C[\{n\}])$ and $\det(C[\{n-1,n\}])$ are nonzero,
and that
$\det(C[\{1\}])$ and $\det(C[\{1,2\}])$
are zero.
Hence, $\ell_1 \ell_2 = \tt SS$.
\epf

Together with Theorems \ref{N...} and \ref{A...},
the next result completes the characterization of
the attainable epr-sequences over $\Z_2$.

\begin{thm}\label{S...}
An epr-sequence starting with $\tt S$ is attainable
by a symmetric matrix over $\Z_2$
if and only if
it has one of the following forms:
\ben
\item $\tt S\OL{S}N\OL{N}$\label{S1};
\item $\tt S\OL{S}A$\label{S2};
\item $\tt S\OL{S}AA$\label{S3};
\item $\tt SS\OL{S}AN$\label{S4};
\item $\tt SASA\OL{SA}$\label{S5};
\item $\tt SASA\OL{SA}A$\label{S6};
\item $\tt SA\OL{SA}N$\label{S7}.
\een
\end{thm}

\bpf
Let $\sigma = \ell_1 \ell_2 \cdots \ell_n$ be
an epr-sequence with $\ell_1 = \tt S$.
Suppose that $\sigma = \epr(B)$, where
$B$ is a symmetric matrix over $\Z_2$.
Since an attainable epr-sequence cannot
end with $\tt S$, $n \geq 2$.
If $n=2$, then $\sigma$ is
$\tt SA$ or $\tt SN$.
Suppose $n \geq 3$.
If $\ell_2 = \tt A$ or $\ell_2 = \tt N$, then
Proposition \ref{SA...} and
the {\tt N}-Even Observation imply that
$\sigma$ is either
$\tt SA\OL{SA}$,
$\tt SA\OL{SA}A$,
$\tt SA\OL{SA}N$, or
$\tt SNN\OL{N}$.
Thus, suppose $\ell_2 = \tt S$.
Hence, by Theorem \ref{ASA anywhere},
$\sigma$ does not contain $\tt ASA$.
Let $k$ be a minimal integer with $2 \leq k \leq n-1$
such that
$\ell_k \ell_{k+1} = \tt SN$ or
$\ell_k \ell_{k+1} = \tt SA$;
in the former case, Observation \ref{obs: S... & A...}
implies that $\sigma = \tt SS\OL{S}N\OL{N}$.
Now consider the latter case, namely
$\ell_k \ell_{k+1} = \tt SA$.
If $n = k+1$, then
$\sigma = \tt SS\OL{S}A$.
Thus, suppose $n \geq k+2$.

We now show that $n = k+2$.
Suppose to the contrary that $n \geq k+3$.
By the {\tt AA} Theorem,
$\ell_{k+2} \neq \tt A$.
If $\ell_{k+2} = \tt N$, then
Observation \ref{obs: S... & A...} implies that
$\sigma$ contains $\tt SANN$,
which is prohibited by Proposition \ref{SAXN};
hence, $\ell_{k+2} = \tt S$, so that
$\ell_k \ell_{k+1}\ell_{k+2} = \tt SAS$.
Then, as $\sigma$ does not contain $\tt ASA$,
and because $\tt SASN$ is prohibited by
Proposition \ref{SAXN},
$\ell_{k+3} = \tt S$, implying that
$\sigma$ contains $\tt ASS$ as a
non-initial subsequence,
a contradiction to Proposition \ref{ASS}.
It follows that $n = k+2$, and
therefore that $\sigma$ is either
$\tt SS\OL{S}AA$ or $\tt SS\OL{S}AN$.

Now, we establish the other direction.
We assume that
the sequence under consideration has order $n$.
The sequence $\tt S\OL{S}N\OL{N}$ is attainable by applying Observation \ref{Obs: Append Ns}(2) to the attainable sequence $\tt A\OL{A}$.
The sequence $\tt S\OL{S}A$ is attainable by
\cite[Observation 2.16]{EPR}.
The attainability of $\tt S\OL{S}AA$ follows
by observing that, by the Inverse Theorem,
the inverse of any symmetric matrix attaining the sequence
$\tt A\OL{S}SA$, which is attainable by
Theorem \ref{A...}, has epr-sequence
$\tt S\OL{S}AA$.
To see that the sequence $\tt SS\OL{S}AN$ is attainable,
observe that the argument above forces
the matrix $G_n$ in Lemma \ref{SS...AN, n even}
to attain this sequence when $n$ is even,
and that it forces the matrix whose existence was established in Lemma \ref{SS...AN, n odd}
to attain this sequence when $n$ is odd.
Finally, the sequences
$\tt SASA\OL{SA}$,
$\tt SASA\OL{SA}A$ and
$\tt SA\OL{SA}N$
are attainable by Proposition \ref{SA...}.
\epf

To conclude, we note that there is no known characterization of the epr-sequences that are attainable
by symmetric matrices over the real field or
any other field besides $\Z_2$.
However, the results of
Theorems \ref{N...}, \ref{A...} and \ref{S...} provide such a characterization for symmetric matrices
over $\Z_2$.


\section*{Acknowledgements}
$\null$
\indent
The author expresses his gratitude to Dr. Leslie Hogben, for introducing him to the topic of pr- and epr-sequences.


\end{document}